\documentstyle[12pt]{article}    
\input{amssym.def}
\input{amssym.tex}
\begin{document}      
\title{$AF$-algebras and topology of  3-manifolds}    

\author{Igor  ~Nikolaev
\footnote{Partially supported 
by NSERC.}}


\date{}    
 \maketitle    
    

\newtheorem{thm}{Theorem}    
\newtheorem{lem}{Lemma}    
\newtheorem{dfn}{Definition}    
\newtheorem{rmk}{Remark}    
\newtheorem{cor}{Corollary}    
\newtheorem{prp}{Proposition}    
\newtheorem{exm}{Example}    
\newtheorem{con}{Conjecture} 
\newtheorem{prb}{Problem}

 \begin{abstract}    
We construct a functor which maps conjugate  pseudo-Anosov automorphisms of a surface to 
the so-called stably isomorphic stationary  $AF$-algebras;  the functor gives new topological 
invariants  of three dimensional manifolds  coming  from  the known  invariants  of 
the $AF$-algebras.  The main invariant is a  triple $(\Lambda, [I], K)$,  where $\Lambda$ 
is an integral order in the real  number field $K$ and $[I]$ the equivalence class of ideals 
in $\Lambda$.

\vspace{7mm}    
    

\vspace{5mm}
{\it AMS (MOS) Subj. Class.:  19K, 46L, 57M.}
\end{abstract}

\section*{Introduction} 
{\bf A. Three dimensional manifolds.}
Let $X$ be an orientable surface of genus $g$. We shall denote by $Mod~(X)$
the mapping class group, i.e. a group of the orientation-preserving
automorphisms of $X$ modulo the normal subgroup of trivial automorphisms. 
Let $\phi\in Mod~(X)$ and consider a  mapping torus 
$M_{\phi}=\left\{ X\times [0,1] ~| \quad (x,0)\mapsto (\phi(x), 1), ~x\in X\right\}$.
The  $M_{\phi}$ is a $3$-dimensional manifold and $M_{\phi}\cong M_{\phi'}$
are homotopy equivalent if and only if $\phi'=\psi\circ\phi\circ\psi^{-1}$ are conjugate
by an automorphism  $\psi\in Mod~(X)$ \cite{Hem1}. Equivalently, $M_{\phi}$ is a surface bundle
over the circle given by the monodromy $\phi$;  such  bundles make by far the most 
interesting, the most complex and the most useful  part of the $3$-dimensional 
topology \cite{Thu1}, p.358.  Each  $\phi\in Mod ~(X)$ is isotopic to an automorphism
$\phi'$, such that either (i) $\phi'$ has a finite order, or (ii) $\phi'$ is a pseudo-Anosov 
automorphism, or else (iii) $\phi'$ is reducible by a system of curves to either type  (i) or (ii) 
\cite{Thu2}. Recall {\it ibid.} that $\phi$ is {\it pseudo-Anosov} if there exist a pair of the stable ${\cal F}_s$
and unstable ${\cal F}_u$ mutually orthogonal measured foliations of the surface $X$,
such that $\phi({\cal F}_s)={1\over\lambda_{\phi}}{\cal F}_s$ 
and $\phi({\cal F}_u)=\lambda_{\phi}{\cal F}_u$, where $\lambda_{\phi}>1$
is a dilatation of $\phi$. The foliations ${\cal F}_s$ and ${\cal F}_u$ are minimal,
uniquely ergodic and describe the automorphism $\phi$ up to a power.  
In the sequel, we shall classify surface bundles $M_{\phi}$,  where $\phi$ is
the pseudo-Anosov  automorphisms of a surface $X$.

\medskip\noindent
{\bf B. The $AF$-algebras (\cite{E}).}
The $C^*$-algebra $A$ is an algebra over complex numbers endowed with the norm 
$a\mapsto ||a||$ and an involution $a\mapsto a^*, a\in A$, such that $A$ is
complete with respect to the norm, and such that $||ab||\le ||a||~||b||$ 
and $||a^*a||=||a||^2$ for every $a,b\in A$. 
Any commutative algebra $A$  is isomorphic to the $C^*$-algebra 
$C_0(X)$ of continuous complex-valued functions on a locally compact Hausdorff space $X$; 
the algebras which are not commutative are deemed as   noncommutative topological
spaces. A {\it stable isomorphism} $A\to A'$ is defined as the (usual) isomorphism
$A\otimes {\cal K}\to A'\otimes {\cal K}$, where ${\cal K}$ is the $C^*$-algebra
of compact operators on a Hilbert space; such an isomorphism corresponds to
 a homeomorphism  between the noncommutative spaces $A$ and $A'$.    
 The matrix algebra $M_n({\Bbb C})$ is an example of noncommutative finite-dimensional
$C^*$-algebra;  a natural generalization  are  approximately finite-dimensional ($AF$-) algebras, 
which  are given by an ascending sequence  
$M_1\buildrel\rm\varphi_1\over\longrightarrow M_2\buildrel\rm\varphi_2\over\longrightarrow\dots$
of finite-dimensional semi-simple $C^*$-algebras $M_i=M_{n_1}({\Bbb C})\oplus\dots\oplus M_{n_k}({\Bbb C})$
and homomorphisms $\varphi_i$  arranged into an infinite graph 
 as follows.  The two sets of vertices $V_{i_1},\dots, V_{i_k}$ and $V_{i_1'},\dots, V_{i_k'}$
are joined by the $b_{rs}$ edges, whenever the summand $M_{i_r}$ contains $b_{rs}$
copies of the summand $M_{i_s'}$ under the embedding $\varphi_i$; as $i\to\infty$,
one gets a {\it Bratteli diagram} of the $AF$-algebra. Such a diagram is defined by
an infinite sequence of {\it incidence matrices}  $B_i=(b_{rs}^{(i)})$.  If the homomorphisms 
$\varphi_1 =\varphi_2=\dots=Const$,  the $AF$-algebra is called {\it stationary}; 
its Bratteli diagram looks like an infinite graph  with the incidence matrix 
$B=(b_{rs})$ repeated over and over again.

\medskip\noindent
{\bf C. The functorialty problem.}
Let $\phi\in Mod~(X)$ be a pseudo-Anosov automorphism;  the problem we seek 
a solution is as follows. Given $\phi$ one assigns  to it  an $AF$-algebra, ${\Bbb A}_{\phi}$,
such  that for  every automorphism $h\in Mod~(X)$  the following diagram commutes:

\begin{picture}(300,110)(-80,-5)
\put(20,70){\vector(0,-1){35}}
\put(130,70){\vector(0,-1){35}}
\put(45,23){\vector(1,0){53}}
\put(45,83){\vector(1,0){53}}
\put(15,20){${\Bbb A}_{\phi}$}
\put(128,20){${\Bbb A}_{\phi'}$}
\put(17,80){$\phi$}
\put(117,80){$\phi'=h\circ\phi\circ h^{-1}$}
\put(60,30){\sf stable}
\put(50,10){\sf isomorphism}
\put(54,90){\sf conjugacy}
\end{picture}

\noindent
In words, if $\phi$ and $\phi'$ are the conjugate pseudo-Anosov automorphisms,  then
the corresponding $AF$-algebras ${\Bbb A}_{\phi}$ and ${\Bbb A}_{\phi'}$ are  stably isomorphic;
the following simple example  indicates,  that the functoriality   problem has a solution.

\medskip\noindent
{\bf D. A model example (case $g=1$).}
Let $T^2$ be  two-dimensional torus and $\phi\in Mod~(T^2)$ an
automorphism given by  non-negative hyperbolic matrix  $A_{\phi}\in SL_2({\Bbb Z})$.
Consider a stationary $AF$-algebra ${\Bbb A}_{\phi}$  given   by the infinite 
Bratteli diagram  with the constant incidence matrix  $B=A_{\phi}$;
it is verified directly,  that $F: \phi\mapsto {\Bbb A}_{\phi}$
is a correctly defined map  on the set of hyperbolic matrices with non-negative entries.
Let us show that if $\phi,\phi'\in Mod~(T^2)$ are  conjugate automorphisms,
then ${\Bbb A}_{\phi}$ is stably isomorphic to ${\Bbb A}_{\phi'}$.
Indeed, let $\phi'=h\circ\phi\circ h^{-1}$ for
an $h\in Mod~(X)$;  then $A_{\phi'}=TA_{\phi}T^{-1}$ for
a matrix $T\in SL_2({\Bbb Z})$ and  $(A_{\phi}')^n=(TA_{\phi}T^{-1})^n=
TA_{\phi}^nT^{-1}$ for any  $n\in {\Bbb N}$.  Recall that 
the $AF$-algebras are stably isomorphic if and only if their Bratteli diagrams contain a 
common block of an arbitrary length;  this claim follows from \cite{E}, Theorem 2.3,
where  the order-isomorphism is replaced by an equivalent condition of Bratteli diagrams having 
the same infinite tail. 
Consider the Bratteli diagrams ${\Bbb A}_{\phi}=\lim_{n\to\infty}A_{\phi}^n$
and  ${\Bbb A}_{\phi'}=\lim_{n\to\infty}TA_{\phi}^nT^{-1}$;  the latter have
a common block of arbitrary length. 
Thus,   ${\Bbb A}_{\phi}\otimes {\cal K}\cong {\Bbb A}_{\phi'}\otimes {\cal K}$,
which  gives a solution to the functoriality problem in the case $g=1$.

\medskip\noindent
{\bf E. The $AF$-algebra ${\Bbb A}_{\phi}$ (case $g\ge 1$).}
Denote by ${\cal F}_{\phi}$ the unstable foliation of a pseudo-Anosov automorphism
$\phi\in Mod~(X)$. For brevity, we assume that ${\cal F}_{\phi}$ is an oriented 
foliation given by the trajectories of a closed $1$-form $\omega\in H^1(X; {\Bbb R})$;
if ${\cal F}_{\phi}$ is  not oriented, the standard double cover construction  
brings it to the oriented case \cite{HuMa1}. 
Let $v^{(i)}=\int_{\gamma_i}\omega$, where $\{\gamma_1,\dots,\gamma_n\}$ be a basis
in the relative homology $H_1(X, Sing~{\cal F}_{\phi}; {\Bbb Z})$, such that
$\theta=(\theta_1,\dots,\theta_{n-1})$ is a vector with  positive coordinates 
$\theta_i=v^{(i+1)} / v^{(1)}$; while  each $\theta_i$ depends on a basis in the 
homology group,  the  ${\Bbb Z}$-module generated by $\theta_i$ does not (lemma \ref{lm1}).  
Consider the Jacobi-Perron continued fraction  of vector  $\theta$ (\cite{BE}):
$$
\left(\matrix{1\cr \theta}\right)=
\lim_{k\to\infty} \left(\matrix{0 & 1\cr I & b_1}\right)\dots
\left(\matrix{0 & 1\cr I & b_k}\right)
\left(\matrix{0\cr {\Bbb I}}\right),
$$
where $b_i=(b^{(i)}_1,\dots, b^{(i)}_{n-1})^T$ is a vector of the non-negative integers,  
$I$ the unit matrix and ${\Bbb I}=(0,\dots, 0, 1)^T$.
We shall denote by ${\Bbb A}_{\phi}$  an (isomorphism class of)  $AF$-algebra given by the 
Bratteli diagram, whose incidence matrices coincide with
$B_k=\left(\small\matrix{0 & 1\cr I & b_k}\right)$ for  $k=1,\dots, \infty$.
Notice  that such a definition coincides with the one given in  the model example
(case $g=1$).  The basic lemma says that ${\Bbb A}_{\phi}$ is a stationary
$AF$-algebra (lemma \ref{lm4}).

\medskip\noindent
{\bf F. The result.}
Let $B$ be the incidence matrix of the stationary algebra ${\Bbb A}_{\phi}$;
let $\lambda_B$ be the Perron-Frobenius eigenvalue of $B$ 
and $(v^{(1)}_B,\dots, v^{(n)}_B)$ the corresponding normalized eigenvector
with  $v^{(i)}_B\in K={\Bbb Q}(\lambda_B)$.
The endomorphism ring of the module   ${\goth m}={\Bbb Z}v^{(1)}_B+\dots+{\Bbb Z}v^{(n)}_B$ 
will be denoted by $\Lambda$. The equivalence class of ideals in the ring $\Lambda$
generated by the ideal ${\goth m}$, we shall write as $[I]$. 
Finally,  let $\Phi$ be a category of all pseudo-Anosov automorphisms
of the surface of genus $g\ge 1$; the morphisms
are the conjugations between the automorphisms. Likewise,
let ${\cal A}$ be a category of all  stationary  $AF$-algebras ${\Bbb A}_{\phi}$,
where  $\phi$ runs the set $\Phi$; the morphisms of ${\cal A}$ are the stable 
isomorphisms among the ${\Bbb A}_{\phi}$'s.   Our main result can be 
expressed as follows.
\begin{thm}\label{thm1}
Let  $F:\Phi\to {\cal A}$  be a map given by the formula $\phi\mapsto {\Bbb A}_{\phi}$.   Then:

\medskip
(i) $F$ is a covariant functor, which   maps the conjugate pseudo-Anosov automorphisms to the stably
isomorphic $AF$-algebras;

\smallskip
(ii) $F^{-1}({\Bbb A}_{\phi})=[\phi]$, where $[\phi]=\{\phi'\in \Phi~|~(\phi')^m=\phi^n, ~m,n\in {\Bbb N}\}$
is the commensurability class of the pseudo-Anosov automorphism $\phi$.  
\end{thm}
\begin{cor}\label{cor1}
The triple  $(\Lambda, [I], K)$ is a homotopy invariant of 
the manifold $M_{\phi}$. 
\end{cor}

\tableofcontents

\section{Preliminaries} 
\subsection{Jacobian of measured foliation} 
Let ${\cal F}$ be a measured foliation on a compact surface
$X$ \cite{Thu2}. For the sake of brevity, we shall  assume
that ${\cal F}$ is an oriented foliation, i.e. given by the 
trajectories of a closed $1$-form $\omega$ on $X$; each non oriented  
foliation is covered by an  oriented one  on a 
surface $\widetilde X$, which is a double cover of $X$ ramified 
at the singular points of the half-integer index of the non-oriented
foliation \cite{HuMa1}.  Let $\{\gamma_1,\dots,\gamma_n\}$ be a
basis in the relative homology group $H_1(X, Sing~{\cal F}; {\Bbb Z})$,
where $Sing~{\cal F}$ is the set of singular points of the foliation ${\cal F}$.
It is well known that $n=2g+m-1$, where $g$ is the genus of $X$ and $m= |Sing~({\cal F})|$. 
The periods of $\omega$ in the above basis we shall write as
$
\lambda_i=\int_{\gamma_i}\omega.
$
By a {\it Jacobian} $Jac~({\cal F})$ of ${\cal F}$, we
understand a ${\Bbb Z}$-module ${\goth m}={\Bbb Z}\lambda_1+\dots+{\Bbb Z}\lambda_n$
regarded as a subset of the real line ${\Bbb R}$. 
\begin{lem}\label{lm1}
The ${\Bbb Z}$-module ${\goth m}$  is independent of the choice of a
basis
\linebreak
in $H_1(X, Sing~{\cal F}; {\Bbb Z})$.
\end{lem}
{\it Proof.} Indeed, let $A=(a_{ij})\in GL_n({\Bbb Z})$ and let
$
\gamma_i'=\sum_{j=1}^na_{ij}\gamma_j
$
be a new basis in $H_1(X, Sing~{\cal F}; {\Bbb Z})$. 
Then using the integration rules:
$\lambda_i'  = \int_{\gamma_i'}\omega = \int_{\sum_{j=1}^na_{ij}\gamma_j}\omega=
 \sum_{j=1}^n\int_{\gamma_j}\omega  = \sum_{j=1}^na_{ij}\lambda_j$. 
To prove that ${\goth m}={\goth m}'$, consider the following equations:
${\goth m}'  = \sum_{i=1}^n{\Bbb Z}\lambda_i' = \sum_{i=1}^n {\Bbb Z} \sum_{j=1}^n a_{ij}\lambda_j=
 \sum_{j=1}^n \left(\sum_{i=1}^n a_{ij}{\Bbb Z}\right)\lambda_j  \subseteq  {\goth m}.
$
Let $A^{-1}=(b_{ij})\in GL_n({\Bbb Z})$ be an inverse to the matrix $A$.
Then $\lambda_i=\sum_{j=1}^nb_{ij}\lambda_j'$ and 
$
{\goth m}  = \sum_{i=1}^n{\Bbb Z}\lambda_i = \sum_{i=1}^n {\Bbb Z} \sum_{j=1}^n b_{ij}\lambda_j'=
\sum_{j=1}^n \left(\sum_{i=1}^n b_{ij}{\Bbb Z}\right)\lambda_j'  \subseteq  {\goth m}'. 
$
Since both ${\goth m}'\subseteq {\goth m}$ and ${\goth m}\subseteq {\goth m}'$, we conclude
that ${\goth m}' = {\goth m}$. Lemma \ref{lm1} follows.
$\square$

\bigskip\noindent
Recall that the measured foliations ${\cal F}$ and ${\cal F}'$ are said to be
{\it topologically conjugate}, if there exists an automorphism $h\in Mod~(X)$,
which sends the leaves of the foliation ${\cal F}$ to the leaves of the foliation ${\cal F}'$.
Note that such an equivalence deals with the topological foliations (i.e. the projective
classes of measured foliations \cite{Thu2}) and does not preserve transversal
measure.  
\begin{lem}\label{lm2}
Let ${\cal F}$ and ${\cal F}'$  be topologically conjugate measured foliations 
on a surface $X$. Then 
$
Jac~({\cal F}')=\mu ~Jac~({\cal F}),
$
where $\mu>0$ is a real number.  
\end{lem}
{\it Proof.} 
Let $h: X\to X$ be an automorphism of the surface $X$. Denote
by $h_*$ its action on $H_1(X, Sing~({\cal F}); {\Bbb Z})$
and by $h^*$ on $H^1(X; {\Bbb R})$ connected  by the formula: 
$
\int_{h_*(\gamma)}\omega=\int_{\gamma}h^*(\omega), ~\forall\gamma\in H_1(X, Sing~({\cal F}); {\Bbb Z}), 
~\forall\omega\in H^1(X; {\Bbb R}).
$
Let $\omega,\omega'\in H^1(X; {\Bbb R})$ be the closed $1$-forms whose
trajectories define the foliations ${\cal F}$ and ${\cal F}'$, respectively.
Since ${\cal F}, {\cal F}'$ are topologically conjugate,
$
\omega'= \mu ~h^*(\omega)
$
for a $\mu>0$. 
Let $Jac~({\cal F})={\Bbb Z}\lambda_1+\dots+{\Bbb Z}\lambda_n$ and 
$Jac~({\cal F}')={\Bbb Z}\lambda_1'+\dots+{\Bbb Z}\lambda_n'$. Then
$
\lambda_i'=\int_{\gamma_i}\omega'=\mu~\int_{\gamma_i}h^*(\omega)=
\mu~\int_{h_*(\gamma_i)}\omega, \qquad 1\le i\le n.
$
By lemma \ref{lm1}, it holds:
$
Jac~({\cal F})=\sum_{i=1}^n{\Bbb Z}\int_{\gamma_i}\omega=
\sum_{i=1}^n{\Bbb Z}\int_{h_*(\gamma_i)}\omega.
$
Therefore
$
Jac~({\cal F}')=\sum_{i=1}^n{\Bbb Z}\int_{\gamma_i}\omega'=
\mu~\sum_{i=1}^n{\Bbb Z}\int_{h_*(\gamma_i)}\omega=\mu~Jac~({\cal F}).
$
Lemma \ref{lm2} follows.
$\square$

\subsection{Functoriality for measured foliations} 
Let ${\cal F}$ be a foliation of surface $X$ endowed with the unique ergodic measure;
suppose that  ${\cal F}$ is   given by the trajectories of a closed $1$-form $\omega\in H^1(X; {\Bbb R})$.
Let $v^{(i)}=\int_{\gamma_i}\omega$, where $\{\gamma_1,\dots,\gamma_n\}$ be a basis
in the relative homology $H_1(X, Sing~{\cal F}_{\phi}; {\Bbb Z})$, such that
$\theta=(\theta_1,\dots,\theta_{n-1})$ is a vector with  the positive coordinates 
$\theta_i=v^{(i+1)} / v^{(1)}$. Consider the  Jacobi-Perron continued
fraction of $\theta$ (\cite{BE}):
$$
\left(\matrix{1\cr \theta}\right)=
\lim_{k\to\infty} \left(\matrix{0 & 1\cr I & b_1}\right)\dots
\left(\matrix{0 & 1\cr I & b_k}\right)
\left(\matrix{0\cr {\Bbb I}}\right),
$$
where $b_i=(b^{(i)}_1,\dots, b^{(i)}_{n-1})^T$ is a vector of the non-negative integers,  
$I$ the unit matrix and ${\Bbb I}=(0,\dots, 0, 1)^T$.
Let  ${\Bbb A}_{\cal F}$ be an (isomorphism class of)  $AF$-algebra given by the 
Bratteli diagram, whose incidence matrices coincide with
$B_k=\left(\small\matrix{0 & 1\cr I & b_k}\right)$ for  all $k=1,\dots, \infty$;
notice  that ${\Bbb A}_{\cal F}$ is correctly defined, since the Jacobi-Perron
fraction of uniquely ergodic measured foliation is convergent \cite{Bau1}.
The following lemma establishes functoriality of the algebras ${\Bbb A}_{\cal F}$
with respect to the topological conjugacy.   
\begin{lem}\label{lm3}
If  ${\cal F}$ and ${\cal F}'$ are topologically conjugate  foliations, then
${\Bbb A}_{\cal F}$ and  ${\Bbb A}_{{\cal F}'}$ are stably isomorphic $AF$-algebras.
\end{lem}
{\it Proof.}
(i) First, let us show that if  ${\goth m}={\Bbb Z}\lambda_1+\dots+{\Bbb Z}\lambda_n$
and   ${\goth m}'={\Bbb Z}\lambda_1'+\dots+{\Bbb Z}\lambda_n'$
are  two ${\Bbb Z}$-modules, such that ${\goth m}'=\mu {\goth m}$ for a $\mu>0$, 
then  the Jacobi-Perron continued fractions of the vectors $\lambda$ and $\lambda'$
coincide except, may be,  a finite number of terms. 
Indeed,  let ${\goth m}={\Bbb Z}\lambda_1+\dots+{\Bbb Z}\lambda_n$ and 
${\goth m}'={\Bbb Z}\lambda_1'+\dots+{\Bbb Z}\lambda_n'$. Since
${\goth m}'=\mu {\goth m}$, where $\mu$ is a positive real,
one gets the following identity of the ${\Bbb Z}$-modules:
$
{\Bbb Z}\lambda_1'+\dots+{\Bbb Z}\lambda_n'={\Bbb Z}(\mu\lambda_1)+\dots+{\Bbb Z}(\mu\lambda_n).
$
One can always assume that $\lambda_i$ and $\lambda_i'$ are positive reals;
 there exists a basis $\{\lambda_1^{''},\dots,\lambda_n^{''}\}$
of the module ${\goth m}'$, such that:
$$
\left\{
\begin{array}{cc}
\lambda'' &= A(\mu\lambda) \nonumber\\
\lambda'' &= A'\lambda',
\end{array}
\right.
$$
where $A,A'\in GL^+_n({\Bbb Z})$ are the matrices, whose entries 
are non-negative integers.  In view of the Proposition 3 of  \cite{Bau1}:
$$
\left\{
\begin{array}{cc}
A &=  \left(\matrix{0 & 1\cr I & b_1}\right)\dots
\left(\matrix{0 & 1\cr I & b_k}\right)\nonumber\\
A' &= \left(\matrix{0 & 1\cr I & b_1'}\right)\dots
\left(\matrix{0 & 1\cr I & b_l'}\right),
\end{array}
\right.
$$
where $b_i, b_i'$ are non-negative integer vectors.
Since the continued fraction for the vectors
$\lambda$ and $\mu\lambda$ coincide for any $\mu>0$ \cite{BE},
we conclude that: 
$$
\left\{
\begin{array}{cc}
\left(\matrix{1\cr \theta}\right)
 &=  \left(\matrix{0 & 1\cr I & b_1}\right)\dots
\left(\matrix{0 & 1\cr I & b_k}\right)
\left(\matrix{0 & 1\cr I & a_1}\right)
\left(\matrix{0 & 1\cr I & a_2}\right)\dots
\left(\matrix{0\cr {\Bbb I}}\right)
\nonumber\\
\left(\matrix{1\cr \theta'}\right)
 &= \left(\matrix{0 & 1\cr I & b_1'}\right)\dots
\left(\matrix{0 & 1\cr I & b_l'}\right)
\left(\matrix{0 & 1\cr I & a_1}\right)
\left(\matrix{0 & 1\cr I & a_2}\right)\dots
\left(\matrix{0\cr {\Bbb I}}\right),
\end{array}
\right.
$$
where 
$$
\left(\matrix{1\cr \theta''}\right)=
\lim_{i\to\infty} \left(\matrix{0 & 1\cr I & a_1}\right)\dots
\left(\matrix{0 & 1\cr I & a_i}\right)
\left(\matrix{0\cr {\Bbb I}}\right). 
$$
In other words, the continued fractions of the vectors $\lambda$ and $\lambda'$
coincide, except a finite number of terms.

\medskip
(ii) By lemma \ref{lm2}  topologically conjugate  foliations ${\cal F}$ and ${\cal F}'$ have
proportional Jacobians, i.e. ${\goth m}'=\mu {\goth m}$
for a $\mu>0$.  Thus, the continued fraction expansion of the basis vectors
of the proportional Jacobians must coincide, except
a finite number of terms; the $AF$-algebras
${\Bbb A}_{\cal F}$ and ${\Bbb A}_{{\cal F}'}$ 
are given by  the Bratteli diagrams, which are identical,
except a finite part of the diagram.  
It is well  known (\cite{E}, Theorem 2.3)  that the $AF$-algebras,  
which have such a  property, are stably isomorphic.      
$\square$

\subsection{Basic lemma} 
There exists a countable family  of measured foliations, which come from the
 pseudo-Anosov automorphisms of surfaces;  we shall  restrict our attention to this class of
foliations.   Let $\phi\in Mod~(X)$ be a pseudo-Anosov automorphism of the surface; then there 
exist a stable ${\cal F}_s$ and unstable ${\cal F}_u$ mutually orthogonal measured foliations on $X$,
such that $\phi({\cal F}_s)={1\over\lambda_{\phi}}{\cal F}_s$ 
and $\phi({\cal F}_u)=\lambda_{\phi}{\cal F}_u$, where $\lambda_{\phi}>1$
is called a dilatation of $\phi$. The foliations ${\cal F}_s,{\cal F}_u$ are minimal,
uniquely ergodic and describe the automorphism $\phi$ up to a power;
we shall understand by ${\cal F}_{\phi}$ the unstable foliation of $\phi$.  Let ${\Bbb A}_{\phi}:={\Bbb A}_{{\cal F}_{\phi}}$
be the $AF$-algebra of the measured foliation ${\cal F}_{\phi}$; the
following lemma describes the basic property of such an algebra  
(to be proved in the next section).
\begin{lem}\label{lm4}
${\Bbb A}_{\phi}$ is stably isomorphic to a stationary $AF$-algebra. 
\end{lem}
Recall that any stationary $AF$-algebra is given by a positive
integer matrix $A$; the similarity class of the matrix corresponds to   the stable
isomorphism class  of the $AF$-algebra ${\Bbb A}_{\phi}$ \cite{E}.

\section{Proofs} 
\subsection{Proof of basic lemma}
Let $\phi\in Mod~(X)$ be a pseudo-Anosov automorphism of the surface $X$;
we proceed by showing, that
invariant foliation ${\cal F}_{\phi}$ is given by form $\omega\in H^1(X; {\Bbb R})$,
which is an eigenvector of the linear map $[\phi]:  H^1(X; {\Bbb R})\to  H^1(X; {\Bbb R})$
induced by $\phi$.
Indeed, let $\lambda_{\phi}$ be a dilatation of $\phi$ and $\Omega$ the corresponding volume element; 
by definition, $\phi(\Omega)=\lambda_{\phi}\Omega$.  Note, that
$\Omega$ is given by restriction of form $\omega$ to a $1$-dimensional 
manifold, transverse to the leaves of ${\cal F}_{\phi}$. The leaves of ${\cal F}_{\phi}$
are fixed by $\phi$ and, therefore, $\phi(\Omega)$ is given by a multiple $\lambda_{\phi}\omega$
of form $\omega$. Since $\omega\in H^1(X; {\Bbb R})$ is a vector, whose
coordinates define ${\cal F}_{\phi}$ up to a scalar, we conclude, that $[\phi](\omega)=\lambda_{\phi}\omega$,
i.e. $\omega$ is an eigenvector of the linear map $[\phi]$.  
Let  $(\lambda_1,\dots,\lambda_n)$ be a basis of the Jacobian of  ${\cal F}_{\phi}$,
such that $\lambda_i>0$. Notice, that $\phi$ acts on $\lambda_i$ as multiplication
by constant $\lambda_{\phi}$; indeed, since $\lambda_i=\int_{\gamma_i}\omega$, we have:
\begin{equation}
\lambda_i'=\int_{\gamma_i}[\phi](\omega)=\int_{\gamma_i}\lambda_{\phi}\omega=\lambda_{\phi}\int_{\gamma_i}\omega=\lambda_{\phi}\lambda_i,
\end{equation}
where $\{\gamma_i\}$ is a basis in $H_1(X, Sing~{\cal F}_{\phi}; {\Bbb Z})$.  Since $\phi$ preserves the leaves
of ${\cal F}_{\phi}$, one concludes that $\lambda_i'\in Jac ~({\cal F}_{\phi})$;
therefore, $\lambda_j'=\sum b_{ij}\lambda_i$ for a non-negative 
integer matrix $B=(b_{ij})$. According to \cite{Bau1}, matrix $B$ can be 
written as a finite product:
\begin{equation}
B=
\left(\matrix{0 & 1\cr I & b_1}\right)\dots
\left(\matrix{0 & 1\cr I & b_p}\right):=B_1\dots B_p,
\end{equation}
where $b_i=(b^{(i)}_1,\dots, b^{(i)}_{n-1})^T$ is a vector of non-negative 
integers and   $I$ the unit matrix.   Let $\lambda=(\lambda_1,\dots,\lambda_n)$.  
Consider a purely periodic Jacobi-Perron continued fraction: 
\begin{equation}
\lim_{i\to\infty} 
\overline{B_1\dots B_p}
\left(\matrix{0\cr {\Bbb I}}\right),
\end{equation}
where  ${\Bbb I}=(0,\dots, 0, 1)^T$;  by a basic property of such fractions, 
it converges to an eigenvector  $\lambda'=(\lambda_1',\dots,\lambda_n')$ of matrix $B_1\dots B_p$
\cite{B},  Ch.3.  But $B_1\dots B_p=B$ and $\lambda$ is an eigenvector of matrix $B$;
therefore, vectors $\lambda$ and $\lambda'$ are collinear.  The collinear vectors  
are known to have the same continued fractions;  thus, we have     
\begin{equation}
\left(\matrix{1\cr \theta}\right)=
\lim_{i\to\infty} 
\overline{B_1\dots B_p}
\left(\matrix{0\cr {\Bbb I}}\right),
\end{equation}
where  $\theta=(\theta_1,\dots,\theta_{n-1})$ and  $\theta_i=\lambda_{i+1}/\lambda_1$. 
Since vector $(1,\theta)$ unfolds into a periodic Jacobi-Perron  continued fraction,
we conclude, that the $AF$-algebra ${\Bbb A}_{\phi}$ is stationary.
Lemma \ref{lm4} is proved.   
$\square$

\subsection{Proof of theorem 1} 
(i) Let us prove the first statement;  we start with the  following lemma. 
\begin{lem}\label{lm6}
Let $\phi$ and $\phi'$ be the conjugate pseudo-Anosov automorphisms
of a surface $X$. Then the invariant foliations ${\cal F}_{\phi}$
and ${\cal F}_{\phi'}$ are topologically conjugate.  
\end{lem}
{\it Proof.} Let $\phi,\phi'\in Mod~(X)$ be conjugate, i.e 
$\phi'=\psi\circ\phi\circ\psi^{-1}$ for an automorphism $\psi\in Mod~(X)$.
Since $\phi$ is the pseudo-Anosov automorphism, there exists  a measured foliation
${\cal F}_{\phi}$,  such that $\phi({\cal F}_{\phi})=\lambda_{\phi}{\cal F}_{\phi}$.
Let us evaluate the automorphism $\phi'$ on the foliation $\psi({\cal F}_{\phi})$:
\begin{eqnarray}\label{eq12}
\phi'(\psi({\cal F}_{\phi}))  &= \psi\circ\phi\circ\psi^{-1}(\psi({\cal F}_{\phi})) &= 
\psi\phi({\cal F}_{\phi})=\nonumber \\
 &= \psi \lambda_{\phi} {\cal F}_{\phi}  &= \lambda_{\phi} (\psi({\cal F}_{\phi})). 
\end{eqnarray}
Thus, ${\cal F}_{\phi'}=\psi({\cal F}_{\phi})$ is the invariant foliation for the 
pseudo-Anosov automorphism $\phi'$ and foliations ${\cal F}_{\phi}$ and ${\cal F}_{\phi'}$
are topologically conjugate.  Note also,  that the pseudo-Anosov automorphisms $\phi$ and $\phi'$ have
 the same dilatation.  
$\square$

\medskip
One can  prove claim (i) of theorem \ref{thm1}; 
let $\phi$ and $\phi'$ be conjugate pseudo-Anosov 
automorphisms. Functor $F$ acts by the formulas $\phi\mapsto {\Bbb A}_{\phi}$
and $\phi'\mapsto {\Bbb A}_{\phi'}$, where ${\Bbb A}_{\phi}$ and  ${\Bbb A}_{\phi'}$
are the $AF$-algebras of the invariant foliations ${\cal F}_{\phi}$ and  ${\cal F}_{\phi'}$. 
In view of lemma \ref{lm6}, foliations ${\cal F}_{\phi}$ and ${\cal F}_{\phi'}$ are
topologically conjugate.  By lemma \ref{lm3}, the $AF$-algebras ${\Bbb A}_{\phi}$
and ${\Bbb A}_{\phi'}$ are stably isomorphic;  claim  (i) is proved.

\bigskip\noindent
(ii) 
Let $\phi\in Mod~(X)$ be a pseudo-Anosov automorphism. Then there exists a unique
measured foliation  ${\cal F}_{\phi}$, such that $\phi({\cal F}_{\phi})=\lambda_{\phi}{\cal F}_{\phi}$,
where $\lambda_{\phi}>1$;  let us evaluate the automorphism
$\phi^2\in Mod~(X)$ on the foliation ${\cal F}_{\phi}$: 
\begin{eqnarray}\label{eq13}
\phi^2({\cal F}_{\phi}) &= \phi (\phi({\cal F}_{\phi})) &= 
\phi(\lambda_{\phi} {\cal F}_{\phi})=\nonumber \\
= \lambda_{\phi} \phi({\cal F}_{\phi}) &= \lambda_{\phi}^2{\cal F}_{\phi}  &= 
\lambda_{\phi^2}{\cal F}_{\phi}, 
\end{eqnarray}
where $\lambda_{\phi^2}:= \lambda_{\phi}^2$. Thus, the foliation ${\cal F}_{\phi}$
is an invariant foliation for the automorphism $\phi^2$ as well;  by induction,
we conclude that ${\cal F}_{\phi}$ is an invariant foliation for the 
automorphism $\phi^n$ for any $n\ge 1$.  Denote by $[\phi]$ the set of all pseudo-Anosov automorphisms
 $\psi$ of $X$,  such that $\psi^m=\phi^n$ for some positive integers
$m$ and $n$.
\begin{lem}\label{lm7}
The  foliation ${\cal F}_{\phi}$ is an invariant foliation for every  
automorphism $\psi\in [\phi]$. 
\end{lem}
{\it Proof.}
Suppose that $\psi\in Mod~(X)$ is a pseudo-Anosov
automorphism, such that $\psi^m=\phi^n$ for some $m\ge 1$ and $\psi\ne\phi$;
then ${\cal F}_{\phi}$ is an invariant foliation for the automorphism 
$\psi$. Indeed,  ${\cal F}_{\phi}$ is an invariant foliation for the 
automorphism $\psi^m$.  If there exists ${\cal F}'\ne {\cal F}_{\phi}$,
such that the foliation  ${\cal F}'$ is an invariant foliation  of $\psi$, then 
the foliation ${\cal F}'$ is an invariant foliation  of the pseudo-Anosov automorphism $\psi^m$. 
Thus, by the uniqueness of invariant foliations, 
${\cal F}'={\cal F}_{\phi}$. 
$\square$

\medskip
In view of lemma \ref{lm7}, one arrives at the following
identities among the $AF$-algebras:
\begin{equation}\label{eq14}
{\Bbb A}_{\phi}={\Bbb A}_{\phi^2}=\dots={\Bbb A}_{\phi^n}=
{\Bbb A}_{\psi^n}=\dots={\Bbb A}_{\psi^2}={\Bbb A}_{\psi}.
\end{equation}
Thus, the functor $F: \Phi\to {\cal A}$ is not  
injective,  since the preimage $F^{-1}$ of the $AF$-algbera
${\Bbb A}_{\phi}$ is  a countable set of pseudo-Anosov
automorphisms $\psi\in [\phi]$  commensurable with the automorphism
$\phi$.

\medskip
Theorem \ref{thm1} is proved.
$\square$

\subsection{Proof of corollary 1}
Lemma \ref{lm4} says that ${\Bbb A}_{\phi}$
is a stationary $AF$-algebra given by a positive  integer matrix $B$. 
By the Perron-Frobenius theory, matrix $B$ has a real eigenvalue $\lambda_B>1$, 
which exceeds the absolute values of all other roots of the characteristic polynomial of $B$;
note that $\lambda_B$ is an algebraic number.  Consider
a real algebraic number field $K={\Bbb Q}(\lambda_B)$ obtained as 
an extension of the field of the rational numbers by  $\lambda_B$. 
Let $(v^{(1)}_B,\dots,v^{(n)}_B)$ be the eigenvector
corresponding to the eigenvalue $\lambda_B$;  one can  normalize the eigenvector so 
that $v^{(i)}_B\in K$. 
Consider the  ${\Bbb Z}$-module
${\goth m}={\Bbb Z}v^{(1)}_B+\dots+{\Bbb Z}v^{(n)}_B$; denote by $\Lambda$ the 
endomorphism ring of ${\goth m}$ and  by $I$ an ideal in the ring $\Lambda$ generated by ${\goth m}$. 
The ring $\Lambda$ is an order in the algebraic number field $K$ and therefore $I$ belongs to an ideal
class in $\Lambda$;  the ideal class of $I$ is denoted by $[I]$. 
The triple $(\Lambda, [I], K)$ is an invariant of the  stable isomorphism class of the  stationary 
$AF$-algebra ${\Bbb A}_{\phi}$  (Handelman \cite{Han1}, \S 5).  
By theorem \ref{thm1}, $(\Lambda, [I], K)$ is an invariant of the conjugacy class of $\phi$ and 
by Hemion \cite{Hem1} of the homotopy class of manifold $M_{\phi}$. 
$\square$

\section{Numerical invariants}
\subsection{Determinant and signature}
One can derive  numerical invariants of the stable isomorphism
classes of stationary $AF$-algebras from the  triple $(\Lambda, [I], K)$;
one such invariant is associated with the trace function on the algebraic number field $K$. 
Recall that $Tr: K\to {\Bbb Q}$ is a linear function on  $K$ such
that $Tr ~(\alpha+\beta)=Tr~(\alpha)+ Tr~(\beta)$ and
$Tr~(a\alpha)=a ~Tr~(\alpha)$ for $\forall\alpha,\beta\in K$ and
$\forall a\in {\Bbb Q}$. 
Let ${\goth m}$ be a full ${\Bbb Z}$-module in the field $K$.   
The trace function defines a symmetric bilinear form
$q(x,y): {\goth m}\times {\goth m}\to {\Bbb Q}$ by 
the formula
$
(x,y)\mapsto Tr~(xy), ~\forall x,y\in {\goth m}.
$
The form 
\begin{equation}\label{eq16}
q(x,y)=\sum_{j=1}^n\sum_{i=1}^ns_{ij}x_iy_j, \qquad\hbox{where} \quad s_{ij}=Tr~(\lambda_i\lambda_j);
\end{equation}
 depends on the basis $\{\lambda_1,\dots,\lambda_n\}$
in the module ${\goth m}$; however,  certain numerical quantities will not depend on the  basis. 
Namely,  consider a symmetric matrix $S=(s_{ij})$ corresponding to the bilinear form $q(x,y)$.
In a new basis  matrix $S$ will take the form $S'=U^TSU$, where ~$U\in GL_n({\Bbb Z})$;
thus  $det~(S')=det~(U^TSU)=det~(U^T) det~(S) det~(U) =det~(S)$.
Therefore, the rational integer
\begin{equation}\label{eq18}
\Delta= det~(Tr~(\lambda_i\lambda_j)),
\end{equation}
does not depend on the
choice of the basis $\{\lambda_1,\dots,\lambda_n\}$ in the
module ${\goth m}$; it is called a {\it determinant} of the bilinear form  $q(x,y)$.
Clearly,   $\Delta$ discerns the modules ${\goth m}$ and  ${\goth m}'$.

\medskip
Finally, recall that the form $q(x,y)$ can be brought by the integer linear
substitutions to the diagonal form:
\begin{equation}\label{eq19}
s_1x_1^2+s_2x_2^2+\dots+s_nx_n^2,
\end{equation}
where $s_i\in {\Bbb Z}-\{0\}$. We let $s_i^+$ be the positive and $s_i^-$ the 
negative entries in the diagonal form.   In view of the law of inertia
for the bilinear forms, the integer number $\Sigma = (\# s_i^+) - (\# s_i^-)$ does not depend
on the choice of basis in the module ${\goth m}$;
it is  called  a {\it signature} of the form. 
Thus, $\Sigma$ discerns the modules ${\goth m}$ and ${\goth m}'$.

\subsection{Numerical invariants of Anosov automorphisms}
Let $K={\Bbb Q}(\sqrt{d})$ be a quadratic extension of the field of rational numbers ${\Bbb Q}$.
Further we suppose that $d$ is a positive square free integer. Let
\begin{equation}
\omega=\cases{{1+\sqrt{d}\over 2} & if $d\equiv 1 ~mod~4$,\cr
               \sqrt{d} & if $d\equiv 2,3 ~mod~4$.}
\end{equation}
\begin{prp}\label{prp1}
Let $f$ be a positive integer. 
Every order  in $K$ has form  
$\Lambda_f={\Bbb Z} +(f\omega){\Bbb Z}$, where 
$f$ is the conductor of $\Lambda_f$.
\end{prp}
{\it Proof.} See \cite{BS} pp. 130-132.
$\square$

\medskip\noindent 
The proposition \ref{prp1} allows to classify the similarity classes of the full modules 
in the field $K$. Indeed, there exists a finite number of ${\goth m}_f^{(1)},\dots,
{\goth m}_f^{(s)}$ of the non-similar full modules in the field  $K$, whose coefficient
ring is the order $\Lambda_f$, cf Theorem 3, Ch 2.7 of \cite{BS}. Thus, proposition \ref{prp1}
gives a finite-to-one classification of the similarity classes of full modules in the field $K$.

Let $\Lambda_f$ be an order  in $K$ with
the  conductor $f$. Under the addition operation, the order  $\Lambda_f$ is a full module, 
which we denote by ${\goth m}_f$.  
Let us evaluate  the invariants $q(x,y)$, $\Delta$  and $\Sigma$ 
on the module  ${\goth m}_f$. To calculate  $(s_{ij})=Tr(\lambda_i\lambda_j)$,  we let 
$\lambda_1=1,\lambda_2=f\omega$. Then:
\begin{eqnarray}
s_{11} &=& 2, \quad s_{12}=a_{21}= f, \quad s_{22}= {1\over 2} f^2(d+1)\quad \hbox{if} \quad d\equiv 1 ~mod~4\nonumber\\
s_{11} &=&  2, \quad s_{12}=s_{21}= 0, \quad s_{22}= 2f^2d \quad \hbox{if} \quad d\equiv 2,3 ~mod~4, 
\end{eqnarray}
and 
\begin{eqnarray}
q(x,y) &=& 2x^2 +2f xy +{1\over 2}f^2(d+1)y^2\quad \hbox{if} \quad d\equiv 1 ~mod~4\nonumber\\
q(x,y) &=&  2x^2+2f^2dy^2\quad  \hbox{if} \quad d\equiv 2,3 ~mod~4.
\end{eqnarray}
Therefore
\begin{equation}
\Delta=\cases{f^2d  & if $d\equiv 1 ~mod~4$,\cr
               4f^2d & if $d\equiv 2,3 ~mod~4$,}
\end{equation}
and $\Sigma=+2$ in the both cases, where $\Sigma=\# (positive) -\# (negative)$ 
entries in the diagonal normal form of $q(x,y)$.

\subsection{Example}
Let us consider a numerical example,  which illustrates 
an advantage  of the above  invariants in comparison to the classical Alexander polynomials.  
Denote by $M_{\phi}$ and $M_{\phi'}$ the  torus bundles given by  the monodromy
\begin{equation}
B=\left(\matrix{5 & 2\cr 2 & 1}\right)\qquad  \hbox{and}
\qquad B'=\left(\matrix{5 & 1\cr 4 & 1}\right), 
\end{equation}
respectively. The Alexander polynomial of three dimensional manifolds  $M_{\phi}$ and $M_{\phi'}$ 
are identical  $\Delta(t)=\Delta'(t)= t^2-6t+1$.  However, the bundles $M_{\phi}$ and $M_{\phi'}$ are 
not homotopy equivalent.

Indeed, the Perron-Frobenius eigenvector of matrix $B$ is  $v_B=(1, \sqrt{2}-1)$
while of the matrix $B'$ is $v_{B'}=(1, 2\sqrt{2}-2)$. The  bilinear forms for the modules 
${\goth m}_B={\Bbb Z}+(\sqrt{2}-1){\Bbb Z}$ and 
${\goth m}_{B'}={\Bbb Z}+(2\sqrt{2}-2){\Bbb Z}$ can be written as
\begin{equation}
q_B(x,y)= 2x^2-4xy+6y^2,\qquad q_{B'}(x,y)=2x^2-8xy+24y^2,
\end{equation}
respectively.  The  modules ${\goth m}_B$ and  ${\goth m}_{B'}$ are not similar in the number field 
$K={\Bbb Q}(\sqrt{2})$, since   their  determinants $\Delta({\goth m}_B)=8$ and 
$\Delta({\goth m}_{B'})=32$ are not equal. Therefore  the matrices $B$ and $B'$ are not similar 
\footnote{The reader may verify this fact using the method of periods, which dates back
to Gauss. First we have to find the fixed points  $Bx=x$ and $B'x=x$, which gives us $x_B=1+\sqrt{2}$ and
$x_{B'}={1+\sqrt{2}\over 2}$, respectively. Then one unfolds the fixed points into a periodic continued fraction,
which gives us $x_B=[2,2,2,\dots]$ and $x_{B'}=[1,4,1,4,\dots]$. Since the period $(2)$ of $x_B$
differs from the period $(1,4)$ of $B'$,  the matrices $B$ and $B'$ are not similar in
$SL(2, {\Bbb Z})$.}
in $SL(2,{\Bbb Z})$.

    

\vskip1cm

\textsc{The Fields Institute for Mathematical Sciences, Toronto, ON, Canada,  
E-mail:} {\sf igor.v.nikolaev@gmail.com}

\smallskip
{\it Current address: 101-315 Holmwood Ave., Ottawa, ON, Canada, K1S 2R2}

\end{document}